\documentclass[11pt]{amsart}

\allowdisplaybreaks[4]

\newtheorem{theorem}{\bf Theorem}[section]

\newtheorem{lemma}{\bf Lemma}[section]
\newtheorem{proposition}{\bf Proposition}[section]

\theoremstyle{definition}
\newtheorem{definition}{\bf Definition}[section]

\theoremstyle{remark}

\numberwithin{equation}{section}

\newcommand{\Z}{\ensuremath{\mathbb{Z}}}
\newcommand{\C}{\ensuremath{\mathbb{C}}}
\newcommand{\R}{\ensuremath{\mathbb{R}}}

\newcommand{\N}{\ensuremath{\mathbb{N}}}

\begin{document}
\title[Characterization of wavelets and MRA wavelets]{Characterization of wavelets and MRA wavelets on local fields of positive characteristic}

\author{Biswaranjan Behera}
\address{(B. Behera) Theoretical Statistics and Mathematics Unit, Indian Statistical Institute, 203 B. T. Road, Kolkata, 700108, India}
\email{biswa@isical.ac.in}

\author{Qaiser Jahan}
\address{(Q. Jahan) Theoretical Statistics and Mathematics Unit, Indian Statistical Institute, 203 B. T. Road, Kolkata, 700108, India}
\email{qaiser\_r@isical.ac.in}
\thanks{Research of the second author is supported by a grant from CSIR, India.}


\date{\today}


\begin{abstract}
We provide a characterization of wavelets on local fields of positive characteristic based on results on affine and quasi affine frames. This result generalizes the characterization of wavelets on Euclidean spaces by means of two basic equations. We also give another characterization of wavelets. Further, all wavelets which are associated with a multiresolution analysis on a such a local field are also characterized.
\end{abstract}

\maketitle


\section{Introduction}

The concept of wavelet is defined and studied extensively in the Euclidean spaces $\R^n$. The characterization of wavelets of $L^2(\R)$ was obtained independently by Wang~\cite{W} and Gripenberg~\cite{Grip} in terms of two basic equations involving the Fourier transform of the wavelet (see also~\cite{HKLS} and~\cite{HW}). This result was generalized to $L^2(\R^n)$ by Frazier, Garrig\'{o}s, Wang, and Weiss~\cite{FGWW} for dilation by 2 and by Calogero~\cite{C} for wavelets associated with a general dilation matrix. Bownik~\cite{B2} provided a new approach to characterizing multiwavelets in $L^2(\R^n)$. This characterization was obtained by using the results about shift invariant systems and quasi-affine systems in~\cite{B1, RS} and~\cite{CSS}.

The notion of multiresolution analysis (MRA) is closely related to wavelets. In fact, it is well known that one can always construct a wavelet from an MRA. But, all wavelets are not obtained in this way. It was proved independently by Gripenberg~\cite{Grip} and Wang~\cite{W} that a wavelet arises from an MRA if and only if its dimension function is 1 a.e. Calogero and Garrig\'{o}s~\cite{CG} gave a characterization of wavelet families arising from biorthogonal MRAs of multiplicity $d$. This result was improved by Bownik and Garrig\'{o}s in~\cite{BG}, where they provided this characterization in terms of the dimension function.

In this article, we will give a characterization of wavelets on local fields of positive characteristic by using the results on affine and quasi-affine frames on such fields obtained in an earlier article~\cite{BJ4}. We will also characterize the wavelets associated with an MRA. For some other aspects of wavelet on such fields, we refer to~\cite{BJ1}, \cite{BJ2} and~\cite{BJ3}.

The local fields are essentially of two type (excluding the connected local fields $\R$ and $\C$), namely, the local fields of characteristic zero and those of positive characteristic. The $p$-adic field $\mathbb{Q}_p$ is a local field of characteristic $0$. Examples of local fields of positive characteristic are the Cantor dyadic group and the Vilenkin $p$-groups. Even though the algebraic structure of local fields of positive characteristic is same as the real number field and its translation set $\{u(k):k\in\N_0\}$ forms a group, it is not true in general that $u(k)+u(l) = u(k+l)$ for nonnegative integers $k$ and $l$ (see section~2 for details). This problem does not show up in the Euclidean case. We have to deal with issues related to this problem separately.

This article is organized as follows. In section~2, we give a brief introduction to the local fields and Fourier analysis on such a field. In section~3, we recall some results on the affine and quasi-affine systems on a local field of positive characteristic and use them to provide a characterization of wavelets. We also give another characterization of wavelets. In section~4, we characterize the MRA wavelets. 


\section{Preliminaries on local fields}

Let $K$ be a field and a topological space. Then $K$ is called a \emph{locally compact field} or a \emph{local field} if both $K^+$ and $K^*$ are locally compact abelian groups, where $K^+$ and $K^*$ denote the additive and multiplicative groups of $K$ respectively.

If $K$ is any field and is endowed with the discrete topology, then $K$ is a local field. Further, if $K$ is connected, then $K$ is either $\R$ or $\mathbb{C}$. If $K$ is not connected, then it is totally disconnected. So by a local field, we mean a field $K$ which is locally compact, nondiscrete and totally disconnected.

We use the notation of the book by Taibleson~\cite{Taib}. Proofs of all the results stated in this section can be found in the books~\cite{Taib} and~\cite{RV}.

Let $K$ be a local field. Since $K^+$ is a locally compact abelian group, we choose a Haar measure $dx$ for $K^+$. If $\alpha\neq 0, \alpha\in K$, then $d(\alpha x)$ is also a Haar measure. Let $d(\alpha x) = |\alpha|dx$. We call $|\alpha|$ the \emph{absolute value} or \emph{valuation} of $\alpha$. We also let $|0| = 0$.

The map $x\rightarrow |x|$ has the following properties:
\begin{itemize}
\item[(a)] $|x| = 0$ if and only if $x=0$;
\item[(b)] $|xy|=|x||y|$ for all $x,y\in K$;
\item[(c)] $|x+y|\leq\max\{|x|, |y| \}$ for all $x,y\in K$.
\end{itemize}
Property (c) is called the \emph{ultrametric inequality}. It follows that 
\begin{equation*}\label{e.max}
|x+y|=\max\{|x|, |y|\} \mbox{if}~ |x|\ne |y|.
\end{equation*}

The set $\mathfrak{D}=\{x\in K : |x|\leq 1\}$ is called the \emph{ring of integers} in $K$. It is the unique maximal compact subring of $K$. Define $\mathfrak{P}=\{x\in K:|x|<1\}$. The set $\mathfrak{P}$ is called the \emph{prime ideal} in $K$. The prime ideal in $K$ is the unique maximal ideal in $\mathfrak{D}$. It is principal and prime.

Since $K$ is totally disconnected, the set of values $|x|$ as $x$ varies over $K$ is a discrete set of the form $\{s^k: k\in\Z\}\cup \{0\}$ for some $s>0$. Hence, there is an element of $\mathfrak{P}$ of maximal absolute value. Let $\mathfrak{p}$ be a fixed element of maximum absolute value in $\mathfrak{P}$. Such an element is called a \emph{prime element} of $K$. Note that as an ideal in $\mathfrak{D},\mathfrak{P}=\left\langle \mathfrak{p}\right\rangle = \mathfrak{p}\mathfrak{D}$.

It can be proved that $\mathfrak{D}$ is compact and open. Hence, $\mathfrak{P}$ is compact and open. Therefore, the residue space $\mathfrak{D}/\mathfrak{P}$ is isomorphic to a finite field $GF(q)$, where $q=p^c$ for some prime $p$ and $c\in\N$. For a proof of this fact we refer to ~\cite{Taib}.

For a measurable subset $E$ of $K$, let $|E|=\int_K\chi_E(x)dx$, where $\chi_E$ is the characteristic function of $E$ and $dx$ is the Haar measure of $K$ normalized so that $|\mathfrak{D}|=1$. Then, it is easy to see that $|\mathfrak{P}|=q^{-1}$ and $|\mathfrak{p}|=q^{-1}$ (see~\cite{Taib}). It follows that if $x\neq 0$, and $x\in K$, then $|x|=q^k$ for some $k\in\mathbb{Z}$.

Let $\mathfrak{D}^*=\mathfrak{D}\setminus\mathfrak{P}=\{x\in K: |x|=1\}$. $\mathfrak{D}^*$ is the group of units in $K^*$. If $x\neq 0$, we can write $x=\mathfrak{p}^k x'$, with $x'\in\mathfrak{D}^*$.

Recall that $\mathfrak{D}/\mathfrak{P}\cong GF(q)$. Let $\mathcal{U}=\{a_i: i=0, 1,\dots, q-1\}$ be any fixed full set of coset representatives of $\mathfrak{P}$ in $\mathfrak{D}$. Let $\mathfrak{P}^k=\mathfrak{p}^k\mathfrak{D}=\{x\in K: |x|\leq q^{-k}\}, k\in\Z$. These are called \emph{fractional ideals}. Each $\mathfrak{P}^k$ is compact and open and is a subgroup of $K^+$ (see~\cite{RV}). 

If $K$ is a local field, then there is a nontrivial, unitary, continuous character $\chi$ on $K^+$. It can be proved that $K^+$ is self dual (see~\cite{Taib}).

Let $\chi$ be a fixed character on $K^+$ that is trivial on $\mathfrak{D}$ but is nontrivial on $\mathfrak{P}^{-1}$. We can find such a character by starting with any nontrivial character and rescaling. We will define such a character for a local field of positive characteristic. For $y\in K$, we define $\chi_y(x)=\chi(yx)$, $x\in K$.

\begin{definition}
If $f\in L^1(K)$, then the Fourier transform of $f$ is the function $\hat f$ defined by
\[
\hat f(\xi)= \int_K f(x)\overline{\chi_{\xi}(x)}~dx.
\]
\end{definition}
Note that
\[
\hat f(\xi)= \int_K f(x)\overline{\chi(\xi x)}~dx = \int_K f(x)\chi(-\xi x)~dx.
\]

Similar to the standard Fourier analysis on the real line, one can prove the following results.
\begin{itemize}
\item[(a)] The map $f\rightarrow \hat f$ is a bounded linear transformation of $L^1(K)$ into $L^{\infty}(K)$, and $\|\hat f\|_{\infty}\leq \|f\|_1$.
\item[(b)] If $f\in L^1(K)$, then $\hat{f}$ is uniformly continuous.
\item[(c)] If $f\in L^1(K)\cap L^2(K)$, then $\|\hat f\|_2=\|f\|_2$.
\end{itemize}

To define the Fourier transform of function in $L^2(K)$, we introduce the functions $\Phi_k$. For $k\in\Z$, let $\Phi_k$ be the characteristic function of $\mathfrak{P}^k$. 

\begin{definition}
For $f\in L^2(K)$, let $f_k=f\Phi_{-k}$ and
\[
\hat{f}(\xi)=\lim\limits_{k\rightarrow\infty}\hat f_k(\xi)
=\lim\limits_{k\rightarrow \infty}\int_{\left|x\right|\leq q^k} f(x)\overline{\chi_{\xi}(x)}~d\xi,
\]
where the limit is taken in $L^2(K)$.
\end{definition}
We have the following theorem (see Theorem 2.3 in~\cite{Taib}).
\begin{theorem}
The fourier transform is unitary on $L^2(K)$.
\end{theorem}

A set of the form $h+\mathfrak{P}^k$ will be called a \textit{sphere} with centre $h$ and radius $q^{-k}$. It follows from the ultrametric inequality that if $S$ and $T$ are two spheres in $K$, then either $S$ and $T$ are disjoint or one sphere contains the other. Also, note that the characteristic function of the sphere $h+\mathfrak{P}^k$ is $\Phi_k(\cdot-h)$ and that $\Phi_k(\cdot-h)$ is constant on cosets of $\mathfrak{P}^k$.

Let $\chi_u$ be any character on $K^+$. Since $\mathfrak{D}$ is a subgroup of $K^+$, the restriction $\chi_{u}|_\mathfrak{D}$ is a character on $\mathfrak{D}$. Also, as characters on $\mathfrak{D}, \chi_u = \chi_v$ if and only if $u-v\in \mathfrak{D}$. That is, $\chi_u=\chi_v$ if $u+\mathfrak{D}=v+\mathfrak{D}$ and $\chi_u\neq \chi_v$ if $(u+\mathfrak{D})\cap (v+\mathfrak{D})=\phi$. Hence, if $\{u(n)\}_{n=0}^{\infty}$ is a complete list of distinct coset representative of $\mathfrak{D}$ in $K^+$, then $\{\chi_{u(n)}\}_{n=0}^{\infty}$ is a list of distinct characters on $\mathfrak{D}$. It is proved in~\cite{Taib} that this list is complete. That is, we have the following proposition.

\begin{proposition}\label{p.com}
Let $\{u(n)\}_{n=0}^{\infty}$ be a complete list of (distinct) coset representatives of $\mathfrak{D}$ in $K^+$. Then $\{\chi_{u(n)}\}_{n=0}^{\infty}$ is a complete list of (distinct) characters on $\mathfrak{D}$. Moreover, it is a complete orthonormal system on $\mathfrak{D}$.
\end{proposition}

Given such a list of characters $\{\chi_{u(n)}\}_{n=0}^{\infty}$, we define the Fourier coefficients of $f\in L^1(\mathfrak{D})$ as
\[
\hat{f}(u(n))=\int_{\mathfrak{D}}f(x)\overline{\chi_{u(n)}(x)}dx.
\]
The series $\sum\limits_{n=0}^{\infty}\hat{f}(u(n))\chi_{u(n)}(x)$ is called the Fourier series of $f$. From the standard $L^2$-theory for compact abelian groups we conclude that the Fourier series of $f$ converges to $f$ in $L^2(\mathfrak{D})$ and Parseval's identity holds:
\[
\int_{\mathfrak{D}}|f(x)|^2dx= \sum\limits_{n=0}^{\infty}|\hat{f}(u(n))|^2.
\]
Also, if $f\in L^1(\mathfrak{D})$ and $\hat f(u(n))=0$ for all $n=0, 1, 2,\dots$, then $f=0$ a. e.

These results hold irrespective of the ordering of the characters. We now proceed to impose a natural order on the sequence $\{u(n)\}_{n=0}^{\infty}$. Note that $\Gamma=\mathfrak{D}/\mathfrak{P}$ is isomorphic to the finite field $GF(q)$ and $GF(q)$ is a $c$-dimensional vector space over the field $GF(p)$. We choose a set $\{1=\epsilon_0, \epsilon_1, \epsilon_2, \dots, \epsilon_{c-1}\}\subset\mathfrak{D}^*$ such that span$\{\epsilon_j\}_{j=0}^{c-1}\cong GF(q)$.
Let $\N_0=\N\cup \{0\}$. For $n\in \N_0$ such that $0\leq n< q$, we have
\[
n=a_0+a_1 p+\dots+a_{c-1} p^{c-1},\quad 0\leq a_k<p, k=0,1,\dots,c-1.
\]
Define
\begin{equation}\label{e.undef1}
u(n)=(a_0+a_1\epsilon_1+\dots+a_{c-1}\epsilon_{c-1})\mathfrak{p}^{-1}.
\end{equation}
Note that $\{u(n):n=0, 1,\dots, q-1\}$ is a complete set of coset representatives of $\mathfrak{D}$ in ${\mathfrak{P}}^{-1}$. Now, for $n\geq 0$, write
\[
n=b_0+b_1q+b_2q^2+\dots+b_sq^s,\quad 0\leq b_k<q, k=0,1,2,\dots,s,
\]
and define
\begin{equation}\label{e.undef2}
u(n)=u(b_0)+u(b_1)\mathfrak{p}^{-1}+\dots+u(b_s)\mathfrak{p}^{-s}.
\end{equation}

This defines $u(n)$ for all $n\in\N_0$. In general, it is not true that $u(m+n)=u(m)+u(n)$. But it follows that
\begin{equation}\label{eq.un}
u(rq^k+s)=u(r)\mathfrak{p}^{-k}+u(s)\quad{\rm if}~r\geq 0, k\geq 0~{\rm and}~0\leq s <q^k.
\end{equation}

In the following proposition we list some properties of $\{u(n)\}$ which will be used later. For a proof, we refer to~\cite{BJ2}.

\begin{proposition}\label{p.un}
For $n\in\N_0$, let $u(n)$ be defined as in \eqref{e.undef1} and \eqref{e.undef2}. Then
\begin{enumerate}
\item[(a)] $u(n)=0$ if and only if $n=0$. If $k\geq 1$, then $|u(n)|=q^k$ if and only if $q^{k-1}\leq n < q^k$;
\item[(b)] $\{u(k): k\in\N_0\}=\{-u(k): k\in\N_0\}$;
\item[(c)] for a fixed $l\in\N_0$, we have $\{u(l)+u(k): k\in\N_0\}=\{u(k): k\in\N_0\}$.
\end{enumerate}
\end{proposition}

For brevity, we will write $\chi_n=\chi_{u(n)}$ for $n\in\N_0$. As mentioned before, $\{\chi_n: n\in\N_0\}$ is a complete set of characters on $\mathfrak{D}$.

Let $K$ be a local field of characteristic $p>0$ and $\epsilon_0, \epsilon_1, \dots, \epsilon_{c-1}$ be as above. We define a character $\chi$ on $K$ as follows~(see~\cite{Zheng}):
\begin{equation}\label{e.chi}
\chi(\epsilon_{\mu}\mathfrak{p}^{-j})=
\left\{
\begin{array}{lll}
\exp(2\pi i/p), & \mu=0~\mbox{and}~j=1,\\
1, & \mu=1,\dots,c-1~\mbox{or}~j\neq 1.
\end{array}
\right.
\end{equation}
Note that $\chi$ is trivial on $\mathfrak{D}$ but nontrivial on $\mathfrak{P}^{-1}$.

In order to be able to define the concepts of multiresolution analysis and wavelet on local fields, we need analogous notions of translation and dilation. Since $\bigcup\limits_{j\in\Z}\mathfrak{p}^{-j}\mathfrak{D}=K$, we can regard $\mathfrak{p}^{-1}$ as the dilation (note that $|\mathfrak{p}^{-1}|=q$) and since $\{u(n): n\in\N_0\}$ is a complete list of distinct coset representatives of $\mathfrak{D}$ in $K$, the set $\{u(n): n\in\N_0\}$ can be treated as the translation set. Note that it follows from Proposition (\ref{p.un}) that the translation set form a subgroup of $K^+$. 

A function $f$ on $K$ will be called \emph{integral-periodic} if
\[
f(x+u(k))=f(x)~\mbox{for all}~k\in\N_0.
\]


\section{The characterization of wavelets}

For $j\in\Z$ and $y\in K$, we define the dilation operator $\delta_j$ and the translation operator $\tau_y$ as follows:
\[
\delta_jf(x) = q^{j/2}f(\mathfrak{p}^{-j}x)  \quad{\rm and} \quad \tau_yf(x) = f(x-y), \quad f\in L^2(K).
\]

\begin{definition}
Let $\Psi = \{\psi^1, \psi^2, \dots, \psi^L\}$ be a finite family of functions in $L^2(K)$. The \emph{affine system} generated by $\Psi$ is the collection 
\[
X(\Psi) = \{\psi^l_{j,k}: 1\leq l \leq L, j\in\Z, k\in\N_0\},
\]
where $\psi^l_{j,k}(x) = q^{j/2}\psi^l(\mathfrak{p}^{-j}x-u(k)) = \delta_j\tau_{u(k)}\psi^l(x)$. The \emph{quasi-affine system} generated by $\Psi$ is 
\[
\tilde{X}(\Psi) = \{\tilde{\psi}^l_{j,k}:1\leq l \leq L, j\in\Z, k\in\N_0\},
\]
where
\begin{equation}\label{f1}
\tilde{\psi}^l_{j,k}(x)=
\left\{
\begin{array}{lll}
\delta_j\tau_{u(k)}\psi^l(x)=q^{j/2}\psi^l(\mathfrak{p}^{-j}x-u(k)), & j\geq 0, k\in\N_0.\\
q^{j/2}\tau_{u(k)}\delta_j\psi^l(x)=q^j\psi^l(\mathfrak{p}^{-j}(x-u(k))), & j<0, k\in\N_0.
\end{array}
\right.
\end{equation}

We say that $\Psi$ is a \emph{set of basic wavelets} of $L^2(K)$ if the affine system $X(\Psi)$ forms an orthonormal basis for $L^2(K)$.
\end{definition}

\begin{definition}
A subset $X$ of $L^2(K)$ is called a \emph{Bessel family} if there exists a constant $B>0$ such that
\begin{eqnarray}
\sum_{\eta\in X}|\langle f, \eta\rangle|^2 \leq B \|f\|^2 \quad ~\mbox{for all}~f\in L^2(K).
\end{eqnarray}
If, in addition, there exists a constant $A>0, A\leq B$ such that
\begin{eqnarray}\label{f2}
A\|f\|^2 \leq \sum_{\eta\in X}|\langle f, \eta\rangle|^2 \leq B \|f\|^2 \quad ~\mbox{for all}~f\in L^2(K),
\end{eqnarray}
then $X$ is called a \emph{frame}. The frame is \emph{tight} if we can choose $A$ and $B$ such that $A=B$.

The affine system $X(\Psi)$ is an \emph{affine frame} if~(\ref{f2}) holds for $X=X(\Psi)$. Similarly, the quasi-affine system $\tilde{X}(\Psi)$ is a \emph{quasi-affine frame} if~(\ref{f2}) holds for $X=\tilde{X}(\Psi)$. 
\end{definition}

Ron and Shen in~\cite{RS} and later Chui, Shi and St\"{o}ckler in~\cite{CSS} have observed the relationship between affine and quasi-affine frames in $\R^n$. In \cite{BJ4}, we have extended their result to the case of local fields of positive characteristic. 

\begin{theorem}\label{thm:affine}
Let $\Psi$ be a finite subset of $L^2(K)$. Then
\begin{itemize}
\item [(a)] $X(\Psi)$ is a Bessel family if and only if $\tilde{X}(\Psi)$ is a Bessel family. Furthermore, their exact upper bounds are equal.
\item [(b)] $X(\Psi)$ is an affine frame if and only if $\tilde{X}(\Psi)$ is a quasi-affine frame. Furthermore, their lower and upper exact bounds are equal.
\end{itemize}
\end{theorem}

\begin{definition}
Given $\{t_i:i\in\N\}\subset\ell^2(\N_0)$, define the operator $H:\ell^2(\N_0)\rightarrow\ell^2(\N)$ by
\[
H(v) = \Big(\langle v, t_i\rangle\Big)_{i\in\N}.
\]
If $H$ is bounded then $\tilde{G} = H^*H:\ell^2(\N_0)\rightarrow\ell^2(\N_0)$ is called the \emph{dual Gramian} of $\{t_i:i\in\N\}$.
\end{definition}

Observe that $\tilde{G}$ is a non negative definite operator on $\ell^2(\N_0)$. Also, note that for $r,s\in\N_0$, we have 
\[
\langle \tilde{G}e_r, e_s \rangle = \langle He_r, He_s\rangle = \sum_{i\in\N}\overline{t_i(r)}t_i(s),
\]
where $\{e_i:i\in\N_0\}$ is the standard basis of $\ell^2(\N_0)$.

The following result characterizes when the system of translates of a given family of functions is a frame in terms of the dual Gramian. The proof is an easy generalization of the corresponding results on the Euclidean cases given in~\cite{RS} and~\cite{B1}.

\begin{theorem}\label{thm:gframe}
Let $\{\varphi_i:i\in\N\}\subset L^2(K)$ and for a.e. $\xi\in\mathfrak{D}$, let $\tilde{G}(\xi)$ denote the dual Gramian of $\{t_i = (\hat{\varphi_i}(\xi+u(k)))_{k\in\N_0}:i\in\N\}\subset\ell^2(\N_0)$. The system of translates $\{T_k\varphi_i:k\in\N_0, i\in\N\}$ is a frame for $L^2(K)$ with constants $A$ and $B$ if and only if $\tilde{G}(\xi)$ is bounded for a.e. $\xi\in\mathfrak{D}$ and
\[
A\|v\|^2\leq \langle\tilde{G}(\xi)v, v\rangle \leq B\|v\|^2 ~\mbox{for}~v\in \ell^2(\N_0) ~\mbox{and for a.e.}~\xi\in\mathfrak{D},
\]
that is, the spectrum of $\tilde{G}(\xi)$ is contained in $[A, B ]$ for a.e. $\xi\in\mathfrak{D}$,
\end{theorem}

We first prove a lemma which gives necessary and sufficient conditions for the orthonormality of an affine system.
\begin{lemma}\label{lem:org}
Suppose that $\Psi = \{\psi^1, \psi^2, \dots \psi^L\}\subseteq L^2(K)$. The affine system $X(\Psi)$ is orthonormal in $L^2(K)$ if and only if 
\begin{equation}\label{e.ortho}
\sum_{k\in\N_0}\hat{\psi}^l(\xi+u(k))\overline{\hat{\psi}^m(\mathfrak{p}^{-j}(\xi+u(k)))} = \delta_{j,0}\delta_{l,m}~\mbox{for a.e.}~\xi\in K, 1\leq l,m\leq L, j\geq 0.
\end{equation}
\end{lemma}

\proof
Using Proposition \ref{p.un}(b) and (c), we observe that 
\[
\langle \psi^l_{j,k}, \psi^{l'}_{j',k'}\rangle = \delta_{l,l'}\delta_{j,j'}\delta_{k,k'}, \quad 1\leq l,l'\leq L,~j,j'\in\Z,~k,k'\in\N_0
\]
is equivalent to 
\[
\langle \psi^l_{j,k}, \psi^{l'}_{0,0}\rangle = \delta_{l,l'}\delta_{j,0}\delta_{k,0} \quad 1\leq l,l'\leq L, j\geq 0, k\in\N_0.
\]
Now, let $1\leq l\leq L$, $j\geq 0$, $k\in\N_0$. Then
\begin{eqnarray*}
\langle \psi^l_{j,k}, \psi^{l'}_{0,0}\rangle
& = & \langle \hat{\psi}^l_{j,k}, \hat{\psi}^{l'}_{0,0}\rangle \\
& = & \int_K q^{-j/2}\hat{\psi}^l(\mathfrak{p}^j\xi)\overline{\chi_k(\mathfrak{p}^j\xi)}\overline{\hat{\psi}^{l'}(\xi)}d\xi \\
& = & \int_K q^{j/2}\hat{\psi}^l(\xi)\overline{\chi_k(\xi)}\overline{\hat{\psi}^{l'}(\mathfrak{p}^{-j}\xi)}d\xi \\
& = & q^{j/2}\int_{\mathfrak{D}}\Big\{\sum_{n\in\N_0}\hat{\psi}^l(\xi+u(n))\overline{\hat{\psi}^{l'}(\mathfrak{p}^{-j}(\xi+u(n)))}\Big\}\overline{\chi_k(\xi)}d\xi.
\end{eqnarray*}
If $\langle\psi^l_{j,k}, \psi^{l'}_{0,0}\rangle=\delta_{l,l'}\delta_{j,0}\delta_{k,0}$ for all $l,l'\in\{1,2,\dots L\}$, $j\geq0$ and $k\in\N_0$, then the $L^1(\mathfrak{D})$ function $F$, where $F(\xi)=\sum\limits_{n\in\N_0}\hat{\psi}^l(\xi+u(n))\overline{\hat{\psi}^{l'}(\mathfrak{p}^{-j}(\xi+u(n)))}$, has the property that its Fourier coefficients are all zero except for the coefficient corresponding to $k=0$, which is $1$ if $j=0$ and $l=l'$. Hence, $F=\delta_{j,0}\delta_{l,l'}$ a.e. Conversely, if $F=\delta_{j,0}\delta_{l,l'}$ a.e, then the same calculation shows that $\langle \psi^l_{j,k}, \psi^{l'}_{0,0}\rangle = \delta_{l,l'}\delta_{j,0}\delta_{k,0}$, since $\{\chi_n:n\in\N_0\}$ is an orthonormal basis of $L^2(\mathfrak{D})$ (see Proposition~\ref{p.com}).
\qed

Define $D_j$ as follows:
\begin{eqnarray*}
D_j=
\left\{
\begin{array}{lll}
\{0,1,\dots q^j-1\}, & j\geq 0,\\
0, & j<0.
\end{array}
\right.
\end{eqnarray*}

Let 
\begin{eqnarray*}
\mathcal{A}
& = & \{\tilde{\psi}^l_{j,d}: 1\leq l\leq L, j\in\Z, d\in D_j\} \\
& = & \{\tilde{\psi}^l_{j,0}: 1\leq l\leq L, j<0\} \cup \{\tilde{\psi}^l_{j,d}: 1\leq l\leq L, j\geq 0, d\in D_j\}.
\end{eqnarray*}

If $j<0$, then $\tau_{u(k)}\tilde{\psi}^l_{j,0}(x)=\tilde{\psi}^l_{j,0}(x-u(k)) = q^j\psi^l(\mathfrak{p}^{-j}(x-u(k))) = \tilde{\psi}^l_{j,k}$. For $j\geq0$, $0\leq d\leq q^j-1$, $k\geq0$, we have
\begin{eqnarray*}
\tau_{u(k)}\tilde{\psi}^l_{j,d}(x)
& = & \tilde{\psi}^l_{j,d}(x-u(k)) = \psi^l_{j,d}(x-u(k)) \\
& = & q^{j/2}\psi^l(\mathfrak{p}^{-j}(x-u(k))-u(d))\\
& = & q^{j/2}\psi^l(\mathfrak{p}^{-j}x-(\mathfrak{p}^{-j}u(k)+u(d)))\\
& = & q^{j/2}\psi^l(\mathfrak{p}^{-j}x-u(kq^j+d))\\
& = & \psi^l_{j,kq^j+d}(x).
\end{eqnarray*}
Since it is true that for each $j\geq0$, every non negative integer $m$ can uniquely be written as $m=kq^j+d$, where $k\in\N_0$, $d\in D_j$, it follows that 
\[
\tilde{X}(\Psi) = \{\tau_{u(k)}\varphi: k\in\N_0, \varphi \in\mathcal{A}\}.
\]

We now define the dual Gramian $\tilde{G}(\xi)$ of the quasi-affine system $\tilde{X}(\Psi)$ at $\xi\in\mathfrak{D}$ to be the dual Gramian of $\{(\hat{\varphi}(\xi+u(k)))_{k\in\N_0}:\varphi\in\mathcal{A}\}\subset\ell^2(\N_0)$. The following lemma will be used later in the computation of $\tilde{G}(\xi)$.

\begin{lemma}\label{lem:chrel}
Let $j\geq 0$. For $p,k\in \N_0$,
\begin{equation}\label{eq.ch}
q^{-j}\sum_{t\in D_j}\chi\Bigl(\bigl(u(p)-u(k)\bigr){\mathfrak p}^ju(t)\Bigr)=
\left\{
\begin{array}{ll}
1, & if \quad p-k\in q^j\N_0.\\
0, & otherwise.
\end{array}
\right.
\end{equation}
\end{lemma}

\proof
The integers $p,k\in \N_0$ can uniquely be written as $p=r+q^jm_1$ and $k=s+q^jm_2$, where $m_1,m_2\in\N_0$ and $0\leq r,s \leq q^j-1$. Using~(\ref{eq.un}), we have $u(p)=u(r)+\mathfrak{p}^{-j}u(m_1)$ and $u(k)=u(s)+\mathfrak{p}^{-j}u(m_2)$. Hence,
\begin{eqnarray*}
\chi\bigl((u(p)-u(k))\mathfrak{p}^ju(t)\bigr) 
& = & \chi\big((u(r)-u(s))\mathfrak{p}^ju(t)+(u(m_1)-u(m_2))u(t)\big) \\
& = & \chi\big((u(r)-u(s))\mathfrak{p}^ju(t)\big),
\end{eqnarray*}
since $\chi(u(k)u(l)) = 1$ for $k,l\in\N_0$ (see Proposition 2 in~\cite{JLJ}). Hence, it is enough to show that if $0\leq r,s\leq q^j-1$, then
\begin{equation}\label{e.chrel}
q^{-j}\sum\limits_{t=0}^{q^j-1}\chi\bigl((u(r)-u(s))\mathfrak{p}^ju(t)\bigr) = \delta_{r,s}.
\end{equation}

If $r=s$ then $u(r)-u(s)=0$, hence both sides of the equation~\eqref{e.chrel} are 1. We now assume that $r\neq s$. Let
\[
r=a_0+a_1q+\dots+a_{j-1}q^{j-1}~{\rm and}~s=b_0+b_1q+\dots+b_{j-1}q^{j-1},
\]
where $0\leq a_m, b_m\leq q-1$ for $m=0, 1, \dots, j-1$. Then from (\ref{e.undef2}), we have
\[
u(r)= u(a_0)+u(a_1)\mathfrak{p}^{-1}+\dots + u(a_{j-1})\mathfrak{p}^{-j+1}
\]
and 
\[
u(s)= u(b_0)+u(b_1)\mathfrak{p}^{-1}+\dots + u(b_{j-1})\mathfrak{p}^{-j+1}.
\]
Similarly, let $t=c_0+c_1q+\dots+c_{j-1}q^{j-1}$, where $0\leq c_n\leq q-1$ for $n=0, 1, \dots, j-1$ so that
\[
u(t)= u(c_0)+u(c_1)\mathfrak{p}^{-1}+\dots + u(c_{j-1})\mathfrak{p}^{-j+1}.
\]

Recall that $\mathfrak{D}/\mathfrak{P}\cong GF(q)\cong {\rm span}\{\epsilon_j\}_{j=0}^{c-1}$. Since $\{u(n)\mathfrak{p}: n=0, 1, \dots, q-1\}$ is a complete set of coset representatives of $\mathfrak{P}$ in $\mathfrak{D}$, for each $n=0, 1, \dots, j-1$, we can write
\[
u(c_n)\mathfrak{p}=\lambda_0^{n}\epsilon_0+\lambda_1^{n}\epsilon_1+\dots+\lambda_{c-1}^{n}\epsilon_{c-1},
\]
where $0\leq\lambda_0^{n}, \lambda_1^{n},\dots,\lambda_{c-1}^{n}\leq p-1$. It can easily be seen that for each $l=0, 1,\dots, c-1$, $\{\epsilon_lu(n)\mathfrak{p}: n=0, 1, \dots, q-1\}$ is also a complete set of coset representatives of $\mathfrak{P}$ in $\mathfrak{D}$. Hence, we have
\[
\epsilon_lu(a_m)\mathfrak{p} = \alpha_0^{m,l}\epsilon_0+\alpha_1^{m,l}\epsilon_1+\dots+\alpha_{c-1}^{m,l}\epsilon_{c-1},\quad l=0, 1,\dots, c-1,
\]
where $0\leq\alpha_0^{m,l}, \alpha_1^{m,l}, \dots, \alpha_{c-1}^{m,l}\leq p-1$. Therefore,
\begin{eqnarray*}
u(r)\mathfrak{p}^ju(t)
& = & \sum_{m=0}^{j-1}\sum_{n=0}^{j-1}\big(u(a_m)\mathfrak{p}u(c_n)\mathfrak{p}\big)\mathfrak{p}^{j-m-n-2} \\
& = & \sum_{m=0}^{j-1}\sum_{n=0}^{j-1}\Big(\sum_{l=0}^{c-1}\lambda_l^n\epsilon_lu(a_m)\mathfrak{p}\Big)\mathfrak{p}^{j-m-n-2} \\
& = & \sum_{m=0}^{j-1}\sum_{n=0}^{j-1}\sum_{l=0}^{c-1}\sum_{k=0}^{c-1}\lambda_l^n\alpha_k^{m,l}\epsilon_k\mathfrak{p}^{j-m-n-2}.
\end{eqnarray*}
By the definition of the character $\chi$ (see~\eqref{e.chi}), we have
\[
\chi(u(r)\mathfrak{p}^ju(t)) = \exp\bigg(\tfrac{2\pi i}{p}\sum_{l=0}^{c-1}(\lambda_l^0\alpha_0^{j-1,l}+\lambda_l^1\alpha_0^{j-2,l}+\dots+\lambda_l^{j-1}\alpha_0^{0,l})\bigg).
\]
Similarly, we can write
\[
\chi(u(s)\mathfrak{p}^ju(t)) = \exp\bigg(\tfrac{2\pi i}{p}\sum_{l=0}^{c-1}(\lambda_l^0\beta_0^{j-1,l}+\lambda_l^1\beta_0^{j-2,l}+\dots+\lambda_l^{j-1}\beta_0^{0,l})\bigg).
\]
where $0\leq\beta_0^{m,l}\leq p-1$ for $m=0, 1,\dots,j-1$ and $l=0, 1,\dots, c-1$.

Observe that as $t$ varies from $0$ to $q^j-1$, the integers $c_0,\dots,c_{j-1}$ all vary from $0$ to $q-1$. Hence, the integers $\lambda_l^n$ vary from $0$ to $p-1$ for $0\leq l\leq c-1$ and $0\leq n\leq j-1$. Therefore,
\begin{eqnarray*}
\lefteqn{
\sum\limits_{t=0}^{q^j-1}\chi\bigl((u(r)-u(s))\mathfrak{p}^ju(t)\bigr) }\\
& = & \sum\limits_{t=0}^{q^j-1}\chi\bigl(u(r)\mathfrak{p}^ju(t)\bigr)\overline{\chi\big(u(s)\mathfrak{p}^ju(t)\bigr)} \\
& = & \Bigg(\sum\limits_{{\lambda_0^0}=0}^{p-1}\exp\Bigl(\tfrac{2\pi i}{p}(\alpha_0^{j-1,0}-\beta_0^{j-1,0})\lambda_0^0\Bigr)\Bigg)\dots
			\Bigg(\sum\limits_{\lambda_0^{j-1}=0}^{p-1}\exp\Bigl(\tfrac{2\pi i}{p}(\alpha_0^{0,0}-\beta_0^{0,0})\lambda_0^{j-1}\Bigr)\Bigg)\\			
&   & \times\Bigg(\sum\limits_{{\lambda_1^0}=0}^{p-1}\exp\Bigl(\tfrac{2\pi i}{p}(\alpha_0^{j-1,1}-\beta_0^{j-1,1})\lambda_1^0\Bigr)\Bigg)\dots
			\Bigg(\sum\limits_{\lambda_1^{j-1}=0}^{p-1}\exp\Bigl(\tfrac{2\pi i}{p}(\alpha_0^{0,1}-\beta_0^{0,1})\lambda_1^{j-1}\Bigr)\Bigg)\\		
&   & \dots \\				
&   & \times\Bigg(\sum\limits_{\lambda_{c-1}^0=0}^{p-1}\exp\Bigl(\tfrac{2\pi i}{p}(\alpha_0^{j-1,c-1}-\beta_0^{j-1,c-1})\lambda_{c-1}^0\Bigr)\Bigg)\dots
			\Bigg(\sum\limits_{\lambda_{c-1}^{j-1}=0}^{p-1}\exp\Bigl(\tfrac{2\pi i}{p}(\alpha_0^{0,c-1}-\beta_0^{0,c-1})\lambda_{c-1}^{j-1}\Bigr)\Bigg).
\end{eqnarray*}
Since $r\neq s$, we have $a_m\neq b_m$ for some $m=0, 1,\dots,j-1$. We claim that there exists some $l\in\{0, 1,\dots,c-1\}$ such that $\alpha_0^{m,l}\neq\beta_0^{m,l}$. If $\alpha_0^{m,l}=\beta_0^{m,l}$ for all $l\in\{0, 1,\dots,c-1\}$, then since $u(a_m)\mathfrak{p}\neq u(b_m)\mathfrak{p}$, we have
\begin{eqnarray*}
GF(q) & = & {\rm span}\{(u(a_m)\mathfrak{p}-u(b_m)\mathfrak{p})\epsilon_l\}_{l=0}^{c-1}\\
      & = & {\rm span}\big\{(\alpha_0^{m,l}-\beta_0^{m,l})\epsilon_0+\dots+(\alpha_{c-1}^{m,l}-\beta_{c-1}^{m,l}) \epsilon_{c-1}\big\}\\
      & \subseteq & {\rm span}\{\epsilon_1,\epsilon_2,\dots,\epsilon_{c-1}\}.
\end{eqnarray*}
This is a contradiction which proves the claim. Now for any $l$ such that $\alpha_0^{m,l}\neq\beta_0^{m,l}$, we observe that 
\[
\sum\limits_{\lambda_l^{j-1-m}=0}^{p-1}\exp\Bigl(\tfrac{2\pi i}{p}(\alpha_0^{m,l}-\beta_0^{m,l})\lambda_l^{j-1-m}\Bigr)
\]
is a factor in the above product. But its value is equal to 
\[
\tfrac{1-\exp\bigl(2\pi i(\alpha_0^{m,l}-\beta_0^{m,l})\bigr)}{1-\exp\bigl(\tfrac{2\pi i}{p}(\alpha_0^{m,l}-\beta_0^{m,l})\bigr)}=0,
\]
since $\alpha_0^{m,l}-\beta_0^{m,l}$ is an integer with absolute value less than $p$. This completes the proof of the lemma.
\qed

For $s\in\N_0\setminus q\N_0$, define the function
\begin{equation}\label{e.ts}
t_s(\xi)=\sum_{l=1}^L\sum_{j=0}^{\infty}\hat{\psi^l}(\mathfrak{p}^{-j}\xi)\overline{\hat{\psi^l}(\mathfrak{p}^{-j}(\xi+u(s)))}.
\end{equation}

In the following lemma we compute the dual Gramian $\tilde{G}(\xi)$ of the quasi-affine system $\tilde{X}(\Psi)$ at $\xi\in\mathfrak{D}$ in terms of the Fourier transforms of functions in $\Psi$.

\begin{lemma}\label{lem:gramian}
Let $\Psi = \{\psi^1, \psi^2, \dots, \psi^L\}\subseteq L^2(K)$ and $\tilde{G}(\xi)$ be the dual Gramian of $\tilde{X}(\Psi)$ at $\xi\in\mathfrak{D}$. Then 
\begin{equation} 
\langle \tilde{G}(\xi)e_k, e_k\rangle = \sum_{l=1}^L\sum_{j\in\Z} |\hat{\psi}^l(\mathfrak{p}^{-j}(\xi+u(k)))|^2\qquad~\mbox{for}~k\in\N_0,
\end{equation}
and
\begin{equation}
\langle \tilde{G}(\xi)e_{k'}, e_k\rangle = t_s(\mathfrak{p}^m\xi+\mathfrak{p}^mk)\qquad~\mbox{for}~k,k'\in\N_0, k\neq k',
\end{equation}
where $m=\max\{j\geq0:k-k'\in q^j\N_0\}$, $s\in\N_0\setminus q\N_0$ is such that $u(s)=\mathfrak{p}^m(u(k')-u(k))$, and $t_s$ is the function defined in~\eqref{e.ts}.
\end{lemma}
 
\proof
For $k,k'\in\N_0$, we have
\begin{eqnarray*}
\langle \tilde{G}(\xi)e_{k'}, e_k\rangle
& = & \sum_{\varphi\in\mathcal{A}}\hat{\varphi}(\xi+u(k))\overline{\hat{\varphi}(\xi+u(k'))}\\
& = & \sum_{l=1}^{L}\sum_{j<0}\hat{\psi}^l(\mathfrak{p}^j(\xi+u(k)))\overline{\hat{\psi}^l(\mathfrak{p}^j(\xi+u(k')))}\\
&   & \qquad +\sum_{l=1}^{L}\sum_{j\geq0}\hat{\psi}^l(\mathfrak{p}^j(\xi+u(k)))\overline{\hat{\psi}^l(\mathfrak{p}^j(\xi+u(k')))}\\
&   & \qquad\times \sum_{d\in D_j}q^{-j}\overline{\chi_{\xi+u(k)}(\mathfrak{p}^ju(d))}\chi_{\xi+u(k')}(\mathfrak{p}^ju(d))).
\end{eqnarray*}
The sum over $D_j$ is equal to
\[
\sum_{d\in D_j}q^{-j}\overline{\chi_{u(k)}(\mathfrak{p}^ju(d))}\chi_{u(k')}(\mathfrak{p}^ju(d))) = \sum_{d\in D_j}q^{-j}\chi\Big(u(d)\mathfrak{p}^j(u(k')-u(k))\Big).
\]
By Lemma \ref{lem:chrel}, this expression is equal to $1$ if $k'-k\in q^j\N_0$ and $0$ otherwise. Therefore, if $k=k'$, then 
\[
\langle \tilde{G}(\xi)e_k, e_k\rangle = \sum_{l=1}^L\sum_{j\in\Z} \mid \hat{\psi}^l(\mathfrak{p}^{-j}(\xi+u(k)))\mid^2 \qquad~\mbox{for}~k\in\N_0,
\]
If $k\neq k'$, let $m=\max\{j\geq0:k-k'\in q^j\N_0\}$. Then, the sum over $D_j$ will contribute $1$ for each $j=0,1\dots m$ and then $0$ from $m+1$ onwards. Thus,
\begin{eqnarray*}
\langle \tilde{G}(\xi)e_{k'}, e_k\rangle
& = & \sum_{l=1}^{L}\sum_{j=-\infty}^m\hat{\psi}^l(\mathfrak{p}^j(\xi+u(k)))\overline{\hat{\psi}^l(\mathfrak{p}^j(\xi+u(k')))}\\
& = & \sum_{l=1}^{L}\sum^{\infty}_{j=-m}\hat{\psi}^l(\mathfrak{p}^{-j}(\xi+u(k)))\overline{\hat{\psi}^l(\mathfrak{p}^{-j}(\xi+u(k')))}\\
& = & \sum_{l=1}^{L}\sum^{\infty}_{j=0}\hat{\psi}^l(\mathfrak{p}^{m-j}(\xi+u(k)))\overline{\hat{\psi}^l(\mathfrak{p}^{m-j}(\xi+u(k')))}\\
& = & \sum_{l=1}^{L}\sum^{\infty}_{j=0}\hat{\psi}^l(\mathfrak{p}^{-j}(\mathfrak{p}^m\xi+\mathfrak{p}^mu(k)))\\
&   & \qquad\times \overline{\hat{\psi}^l(\mathfrak{p}^{-j}(\mathfrak{p}^m\xi+\mathfrak{p}^mu(k)+\mathfrak{p}^m(u(k')-u(k)))}.
\end{eqnarray*}

Let $k'=c_0+c_1q+\dots +c_Jq^J$ and $k=d_0+d_1q+\dots +d_Jq^J$. Since $k'-k\in q^m\N_0$ but $k'-k\notin q^{m+1}\N_0$, we have $c_i=d_i$ for all $i=1,2,\dots,m-1$ and $c_m \neq d_m$. Hence,
\begin{eqnarray*}
\mathfrak{p}^m(u(k')-u(k)) 
& = & u(c_m+c_{m+1}q+\dots + c_Jq^{J-m})-u(d_m+d_{m+1}q+\dots + d_Jq^{J-m}) \\
& = & u(s),
\end{eqnarray*}
for some $s\in\N_0$, by Proposition~\ref{p.un}(b) and (c). Note that $s\notin q\N_0$, since, otherwise $u(s)=\mathfrak{p}^{-1}u(n)$ for some $n\in\N_0$. This will imply that $c_m=d_m$, which is false. Therefore,
\[
\langle\tilde{G}(\xi)e_{k'}, e_k\rangle = t_s(\mathfrak{p}^m\xi+\mathfrak{p}^mu(k)),
\]
where $s\in\N_0\setminus q\N_0$ is as defined above. This completes the proof of the lemma.
\qed

In the following theorem, we provide ncessary and sufficient conditions for the affine system $X(\Psi)$ to be a tight frame in $L^2(K)$. As a consequence, we get a characterization of wavelets.

\begin{theorem}\label{thm:chr}
Suppose $\Psi = \{\psi^1, \psi^2, \dots, \psi^L\}\subseteq L^2(K)$. The affine system $X(\Psi)$ is a tight frame with constant $1$ for $L^2(K)$, i.e,
\[
\|f\|_2^2 = \sum_{l=1}^L\sum_{j\in\Z}\sum_{k\in\N_0}|\langle f, \psi_{l,j,k}\rangle|^2 \quad~\mbox{for all}~f\in L^2(K).
\]
if and only if the functions $\psi^1, \psi^2, \dots, \psi^L$ satisfy the following two conditions:
\begin{equation}\label{f4}
\sum_{l=1}^L\sum_{j\in\Z}|\hat{\psi^l}(\mathfrak{p}^{-j}\xi)|^2 = 1 \quad~\mbox{for a.e.}~\xi\in K,
\end{equation}
and
\begin{equation}\label{f5}
t_m(\xi)=0\quad~\mbox{for a.e.}~\xi\in K~\mbox{and for all}~m\in\N_0\setminus q\N_0.
\end{equation}
In particular, $\Psi$ is a set of basic wavelets of $L^2(K)$ if and only if $\|\psi^l\|_2 =1$ for $l=1,2,\dots, L$ and ~(\ref{f4}) and ~(\ref{f5}) hold.
\end{theorem}

\proof
It follows from Theorem~\ref{thm:affine} that $X(\Psi)$ is a tight frame with constant $1$ if and only if $\tilde{X}(\Psi)$ is a tight frame with constant $1$. By Theorem \ref{thm:gframe}, this is equivalent to the spectrum of $\tilde{G}(\xi)$ consisting of a single point $1$, i.e., $\tilde{G}(\xi)$ is identity on $\ell^2(\N_0)$ for a.e. $\xi\in\mathfrak{D}$. By Lemma~\ref{lem:gramian}, this is equivalent to the fact that equations (\ref{f4}) and~(\ref{f5}) hold. The second assertion follows since a tight frame $X(\Psi)$ is an orthonormal basis if and only if $\|\psi^l\|_2=1$ for $l=1,2,\dots, L$ (see Theorem~1.8, section~7.1 in~\cite{HW}).
\qed

The following theorem, initially proved by Bownik \cite{B2} for $\R^n$ with an integer dilation matrix, gives a new characterization of tight wavelet frames with constant 1. We extend this result to the case of local fields of positive characteristic.
\begin{theorem}\label{thm:gen}
Suppose $\Psi = \{\psi^1, \psi^2, \dots, \psi^L\}\subseteq L^2(K)$. Assume that $X(\Psi)$ is a Bessel family with constant $1$. Then the following are equivalent:
\begin{itemize}
\item [(a)] $X(\Psi)$ is a tight frame with constant 1.
\item [(b)] $\Psi$ satisfies $(\ref{f4})$.
\item [(c)] $\Psi$ satisfies
\begin{equation}\label{f8}
\sum_{l=1}^L\int_K \frac{|\hat{\psi^l}(\xi)|^2}{|\xi|}d\xi = \frac{q-1}{q}.
\end{equation}
\end{itemize}
\end{theorem}

\proof
It is obvious from Theorem \ref{thm:chr} that (a) $\Rightarrow$ (b). To show (b) implies (c), assume that~\eqref{f4} holds. Then, since $\{\mathfrak{p}^j\mathfrak{D}^*:j\in\Z\}$ is a partition of $K$, we have
\begin{eqnarray*}
\sum_{l=1}^L\int_K|\hat{\psi^l}(\xi)|^2\frac{d\xi}{|\xi|}
& = & \sum_{l=1}^L\sum_{j\in\Z}\int_{\mathfrak{p}^j\mathfrak{D}*}|\hat{\psi^l}(\xi)|^2\frac{d\xi}{|\xi|}\\
& = & \sum_{l=1}^L\sum_{j\in\Z}\int_{\mathfrak{D}*}|\hat{\psi^l}(\mathfrak{p}^{-j}\xi)|^2\frac{q^jd\xi}{|\mathfrak{p}^{-j}\xi|}\\
& = & \int_{\mathfrak{D}*}\Big(\sum_{l=1}^L\sum_{j\in\Z}|\hat{\psi^l}(\mathfrak{p}^{-j}\xi)|^2\Big)\frac{d\xi}{|\xi|}\\
& = & \int_{\mathfrak{D}*}\frac{d\xi}{|\xi|} = |\mathfrak{D}^*|\\
& = & \frac{q-1}{q}.
\end{eqnarray*}
  
To prove (c) $\Rightarrow$ (a), we assume that~(\ref{f8}) holds. Since $X(\Psi)$ is a Bessel family with constant $1$, so is $\tilde{X}(\Psi)$, by Theorem~\ref{thm:affine}(a). Let $\tilde{G}(\xi)$ be the dual Gramian of $\tilde{X}(\Psi)$ at $\xi\in\mathfrak{D}$. By Theorem \ref{thm:gframe}, we have $\|\tilde{G}(\xi)\| \leq 1$ for a.e. $\xi\in\mathfrak{D}$. In particular, $\|\tilde{G}(\xi)e_k\|\leq 1$. Hence,
\begin{eqnarray}\label{f9}
1 \geq \|\tilde{G}(\xi)e_k\|^2 = \sum_{p\in \N_0} |\langle\tilde{G}(\xi)e_k, e_p\rangle|^2 =  |\langle\tilde{G}(\xi)e_k, e_k\rangle|^2 + \sum_{p\in \N_0, p\neq k} |\langle\tilde{G}(\xi)e_k, e_p\rangle |^2 \nonumber.
\end{eqnarray}
By Lemma (\ref{lem:gramian}), we have
\[
\sum_{l=1}^L\sum_{j\in\Z}|\hat{\psi}^l(\mathfrak{p}^{-j}(\xi+u(k)))|^2 = \langle \tilde{G}(\xi)e_k, e_k\rangle \leq 1 ~\mbox{for}~k\in\N_0,~\mbox{a.e.}~\xi\in\mathfrak{D}.
\]
Hence, 
\[
\frac{q-1}{q} = \sum_{l=1}^L\int_K \frac{|\hat{\psi^l}(\xi)|^2}{|\xi|}d\xi = \int_{\mathfrak{D}*}\Big(\sum_{l=1}^L\sum_{j\in\Z}|\hat{\psi^l}(\mathfrak{p}^{-j}\xi)|^2\Big)\frac{d\xi}{|\xi|} \leq \int_{\mathfrak{D}*}\frac{d\xi}{|\xi|} = \frac{q-1}{q}.
\]
From this it follows that $\sum_{l=1}^L\sum_{j\in\Z}|\hat{\psi}^l(\mathfrak{p}^{-j}\xi)|^2=1$ for a.e. $\xi\in\mathfrak{D}^*$ and hence for a.e. $\xi\in K$, i.e., equation~(\ref{f4}) holds. By Lemma~\ref{lem:gramian} and (\ref{f4}), $|\langle \tilde{G}(\xi)e_k, e_k\rangle|^2 = 1$ for all $k\in\N_0$. Hence by~(\ref{f9}), it follows that $\langle\tilde{G}(\xi)e_k, e_{k'}\rangle = 0$ for $k\neq k'$ so that $\tilde{G}(\xi)$ is the identity operator on $\ell^2(\N_0)$. Hence, by Theorem~\ref{thm:gframe}, $\tilde{X}(\Psi)$ is a tight frame with constant $1$. Therefore, $X(\Psi)$ is also a tight frame with constant $1$, by Theorem~\ref{thm:affine}.
\qed

As a consequence of the above theorem, we get a new characterization of wavelets.
\begin{theorem}\label{thm:nchr}
Suppose $\Psi = \{\psi^1, \psi^2, \dots, \psi^L\}\subseteq L^2(K)$. Then the following are equivalent:
\begin{itemize}
\item [(a)] $\Psi$ is a set of basic wavelets of $L^2(K)$.
\item [(b)] $\Psi$ satisfies~\eqref{e.ortho} and~\eqref{f4}.
\item [(c)] $\Psi$ satisfies~\eqref{e.ortho} and~\eqref{f8}.
\end{itemize}
\end{theorem}

\proof
It follows from Theorem~\ref{thm:gen} and Lemma~\ref{lem:gramian} that (a) $\Rightarrow$ (b) $\Rightarrow$ (c). We now prove that (c) implies (a). Assume that $\Psi$ satisfies~(\ref{e.ortho}) and~(\ref{f8}). The equation~(\ref{e.ortho}) implies that $X(\Psi)$ is an orthonormal system, hence it is a Bessel family with constant $1$. By Theorem~\ref{thm:gen} and~(\ref{f8}), $X(\Psi)$ is a tight frame with constant $1$. Since each $\psi^l$ has $L^2$ norm $1$, it follows that $X(\Psi)$ is an orthonormal basis for $L^2(K)$. That is, $\Psi$ is a set of basic wavelets of $L^2(K)$.
\qed


\section{The Characterization of MRA Wavelets}
Similar to $\R^n$, wavelets can be constructed from a multiresolution analysis (MRA). We define an MRA on local fields of positive characteristic as follows (see~\cite{JLJ}).
\begin{definition} \label{MRA}
Let $K$ be a local field of characteristic $p>0$, $\mathfrak{p}$ be a prime element of $K$ and $u(n)\in K$ for $n\in\N_0$ be as defined in~(\ref{e.undef1}) and~(\ref{e.undef2}). A multiresolution analysis (MRA) of $L^2(K)$ is a sequence $\{V_j: j\in\mathbb{Z}\}$ of closed subspaces of $L^2(K)$ satisfying the following properties:
\begin{enumerate}
\item[(a)] $V_j\subset V_{j+1}$ for all $j\in\Z$;
\item[(b)] $\bigcup\limits_{j\in\Z}V_j$ is dense in $L^2(K)$; 
\item[(c)] $\bigcap\limits_{j\in\Z}V_j = \{0\}$;
\item[(d)] $f\in V_j$ if and only if $f(\mathfrak{p}^{-1}\cdot)\in V_{j+1}$ for all $j\in\Z$;
\item[(e)] there is a function $\varphi\in V_0$, called the \emph{scaling function}, such that $\{\varphi(\cdot-u(k)): k\in\N_0\}$ forms an orthonormal basis for $V_0$.
\end{enumerate}
\end{definition}

Let $\Psi=\{\psi^1, \psi^2, \dots, \psi^L\}$ be a set of basic wavelets of $L^2(K)$. We define the spaces $W_j$, $j\in\Z$, by $W_j=\overline{{\rm span}}\{\psi^l_{j, k}:1\leq l\leq L, k\in\N_0\}$. We also define $V_j=\bigoplus\limits_{m<j}W_m$, $j\in\Z$. Then it follows that $\{V_j:j\in\Z\}$ satisfies the properties (a)-(d) in the definition of an MRA. Hence, $\{V_j:j\in\Z\}$ will form an MRA of $L^2(K)$ if we can find a function $\varphi\in L^2(K)$ such that the system $\{\varphi(\cdot-u(k)): k\in\N_0\}$ is an orthonormal basis for $V_0$. In this case, we say that $\Psi$ \emph{is associated with an MRA}, or simply that $\Psi$ is an \emph{MRA-wavelet}.

Now suppose that $\{\psi^1, \psi^2, \dots, \psi^{q-1}\}$ is a set of basic wavelets for $L^2(K)$ associated with an MRA $\{V_j: j\in\Z\}$. Let $\varphi\in L^2(K)$ be the corresponding scaling function. In Theorem~5.1 of~\cite{BJ2}, we have characterized the scaling functions for MRAs of $L^2(K)$. In view of this theorem, we have 
\begin{equation}\label{ons}
\sum\limits_{k\in\N_0}\arrowvert\hat{\varphi}(\xi+u(k))\arrowvert ^2 = 1 \quad{\rm for~a.e.}~\xi\in\mathfrak{D},
\end{equation}
\begin{equation}\label{den}
\lim\limits_{j\rightarrow\infty}\arrowvert\hat{\varphi}(\mathfrak{p}^j\xi)\arrowvert = 1 \quad{\rm for~a.e.}~\xi\in K,
\end{equation}
and
\begin{equation}\label{eqn1}
\hat{\varphi}(\xi) = m_0(\mathfrak{p}\xi)\hat{\varphi}(\mathfrak{p}\xi) \quad{\rm for~a.e.}~\xi\in K,
\end{equation}
where $m_0$ is an integral-periodic function in $L^2(\mathfrak{D})$. Also, since $\{\psi^1, \psi^2, \dots, \psi^{q-1}\}$ are the wavelets associated with an MRA corresponding to the scaling function $\varphi$, there exist integral-periodic functions $m_l$, $1\leq l\leq q-1$, such that the matrix
\[
M(\xi)=\Bigl[m_{l_1}(\mathfrak{p}\xi+\mathfrak{p}u(l_2))\Bigr]_{l_1,l_2=0}^{q-1}
\]
is unitary for a.e. $\xi\in\mathfrak{D}$ (see \cite{JLJ}) and
\[
\hat{\psi^l}(\xi) = m_l(\mathfrak{p}\xi)\hat{\varphi}(\mathfrak{p}\xi).
\]
Hence, by $~(\ref{eqn1})$, we have
\[
|\hat{\varphi}(\xi)|^2 + \sum_{l=1}^{q-1}|\hat{\psi^l}(\xi)|^2 = |m_0(\mathfrak{p}\xi)\hat{\varphi}(\mathfrak{p}\xi)|^2 + \sum_{l=1}^{q-1}|m_l(\mathfrak{p}\xi)\hat{\varphi}(\mathfrak{p}\xi)|^2 = |\hat{\varphi}(\mathfrak{p}\xi)|^2\Big(\sum_{l=0}^{q-1}|m_l(\mathfrak{p}\xi)|^2\Big).
\]
Since $M(\xi)$ is unitary, we have
\[
 |\hat{\varphi}(\xi)|^2 + \sum_{l=1}^{q-1}|\hat{\psi^l}(\xi)|^2 = |\hat{\varphi}(\mathfrak{p}\xi)|^2.
\]
This equality holds for a.e. $\xi\in K$. Hence, we have
\[
|\hat{\varphi}(\xi)|^2 = |\hat{\varphi}(\mathfrak{p}^{-1}\xi)|^2 + \sum_{l=1}^{q-1}|\hat{\psi^l}(\mathfrak{p}^{-1}\xi)|^2.
\]
Iterating, we get, for any integer $N\geq 1$,
\[
|\hat{\varphi}(\xi)|^2 = |\hat{\varphi}(\mathfrak{p}^{-N}\xi)|^2 + \sum_{l=1}^{q-1}\sum_{j=1}^N|\hat{\psi^l}(\mathfrak{p}^{-j}\xi)|^2.
\]
Since $|\hat{\varphi}(\xi)|\leq 1$, the sequence $\{\sum_{j=1}^N\sum_{l=1}^{q-1}|\hat{\psi}^l(\mathfrak{p}^{-j}\xi)|^2:N\geq 1\}$ of real numbers is increasing and is bounded by $1$, hence it converges. Therefore, 
$\lim\limits_{N\rightarrow\infty}|\hat{\varphi}(\mathfrak{p}^{-N}\xi)|^2$ also exists. Now,
\[
\int_K|\hat{\varphi}(\mathfrak{p}^{-N}\xi)|^2d\xi = q^{-N}\int_K|\hat{\varphi}(\xi)|^2d\xi \rightarrow 0 ~\mbox{as}~ N\rightarrow \infty.
\]
Hence, by Fatou's lemma,
\[
\int_K\lim\limits_{N\rightarrow\infty}|\hat{\varphi}(\mathfrak{p}^{-N}\xi)|^2d\xi \leq  \lim\limits_{N\rightarrow\infty}\int_K|\hat{\varphi}(\mathfrak{p}^{-N}\xi)|^2d\xi = 0.
\]
This shows that $\lim\limits_{N\rightarrow\infty}|\hat{\varphi}(\mathfrak{p}^{-N}\xi)|^2 = 0$. Hence, we get
\[
|\hat{\varphi}(\xi)|^2 =  \sum_{l=1}^{q-1}\sum_{j=1}^{\infty}|\hat{\psi^l}(\mathfrak{p}^{-j}\xi)|^2.
\]

Since $\{\varphi(\cdot-u(k)):k\in\N_0\}$ is an orthonormal system, we get for a.e. $\xi\in K$,
\begin{eqnarray}\label{f14}
1=\sum_{k\in\N_0}|\hat{\varphi}(\xi+u(k))|^2=\sum_{l=1}^{q-1}\sum_{j=1}^{\infty}\sum_{k\in\N_0}|\hat{\psi^l}(\mathfrak{p}^{-j}(\xi+u(k)))|^2.
\end{eqnarray}

\begin{definition}
Suppose $\Psi=\{\psi^1,\psi^2,\dots,\psi^L\}\subseteq L^2(K)$ is a set of basic wavelets for $L^2(K)$. The \emph{dimension function} of $\Psi$ is defined as
\[
D_{\Psi}(\xi) = \sum_{l=1}^L\sum_{j=1}^{\infty}\sum_{k\in\N_0}|\hat{\psi^l}(\mathfrak{p}^{-j}(\xi+u(k)))|^2 ~\mbox{a.e.}~\xi\in K.
\] 
\end{definition}

Observe that if $\psi^1,\psi^2,\dots,\psi^L\in L^2(K)$, then 
\begin{equation}\label{f10}
\int_{\mathfrak{D}}\sum_{j=1}^{\infty}\sum_{k\in\N_0}|\hat{\psi^l}(\mathfrak{p}^{-j}(\xi+u(k)))|^2d\xi = \sum_{j=1}^{\infty} q^{-j}\int_K |\hat{\psi^l}(\xi)|^2d\xi < \infty.
\end{equation}
Hence, $D_\Psi$ is well-defined for a.e. $\xi\in K$. In particular, $\sum_{k\in\N_0}|\hat{\psi^l}(\mathfrak{p}^{-j}(\xi+u(k)))|^2 <\infty$ for a.e. $\xi\in K$. Thus, for all $j\geq 1$, $1\leq l\leq L$, and a.e. $\xi\in K$, we can define the vector $\omega_j^l(\xi)$ in $\ell^2(\N_0)$, where
\[
\omega^l_j(\xi) = \{\hat{\psi^l}(\mathfrak{p}^{-j}(\xi+u(k))):k\in\N_0\}.
\]
Note that we can also write $D_\Psi$ as
\begin{equation*}\label{f11}
D_{\Psi}(\xi) = \sum_{l=1}^L\sum_{j=1}^{\infty}\|\omega^l_j(\xi)\|^2_{\ell^2(\N_0)}.
\end{equation*}

We have thus proved that if $\Psi=\{\psi^1,\psi^2,\dots,\psi^{q-1}\}$ is a set of basic wavelets associated with an MRA of $L^2(K)$, then it is necessary that $D_{\Psi}=1$ a.e. Our aim is to show that this condition is also sufficient. We will show that if $\Psi=\{\psi^1,\psi^2,\dots,\psi^{q-1}\}$ is a set of basic wavelets of $L^2(K)$ and $D_{\Psi}=1$ a.e., then $\Psi$ is an MRA-wavelet. To prove this we need the following lemma.

\begin{lemma}\label{lem:rel}
For all $j\geq 1$, $l=1,2,\dots,q-1$, and a.e. $\xi\in K$, we have
\begin{eqnarray}\label{omega}
\omega^l_j(\xi) = \sum_{h=1}^{q-1}\sum_{i=1}^{\infty}\langle \omega^l_j(\xi), \omega^h_i(\xi)\rangle \omega^h_i(\xi).
\end{eqnarray}
\end{lemma}

\proof
The series appearing in the lemma converges absolutely by ~(\ref{f11}) for a.e. $\xi\in K$. We first show that
\begin{eqnarray}\label{f12}
\hat{\psi}^l(\mathfrak{p}^{-j}\xi) = \sum_{h=1}^{q-1}\sum_{i=1}^{\infty}\sum_{k\in\N_0}\hat{\psi^l}(\mathfrak{p}^{-j}(\xi+u(k)))\overline{\hat{\psi}^h(\mathfrak{p}^{-i}(\xi+u(k)))}\hat{\psi}^h(\mathfrak{p}^{-i}\xi)
\end{eqnarray}
Let us denote the series on the right of (\ref{f12}) by $G^l_j(\xi)$. Then by using Lemma \ref{lem:org} and equation (\ref{f5}), we have  
\begin{eqnarray*}
G^l_j(\xi)
& = & \sum_{k\in\N_0}\hat{\psi^l}(\mathfrak{p}^{-j}(\xi+u(k)))\sum_{h=1}^{q-1}\sum_{i=1}^{\infty}\overline{\hat{\psi}^h(\mathfrak{p}^{-i}(\xi+u(k)))}\hat{\psi}^h(\mathfrak{p}^{-i}\xi)\\
& = & \sum_{k\in\N_0}\hat{\psi^l}(\mathfrak{p}^{-j}(\xi+u(k)))\Big(t_k(\xi) - \sum_{h=1}^{q-1}\overline{\hat{\psi}^h(\xi+u(k))}\hat{\psi}^h(\xi)\Big)\\
& = & \sum_{k\in q\N_0}\hat{\psi^l}(\mathfrak{p}^{-j}(\xi+u(k)))t_k(\xi)\\
& = & \sum_{h=1}^{q-1}\sum_{i=0}^{\infty}\sum_{k\in\N_0}\hat{\psi^l}(\mathfrak{p}^{-j}(\xi+u(qk)))\overline{\hat{\psi}^h(\mathfrak{p}^{-i}(\xi+u(qk)))}\hat{\psi}^h(\mathfrak{p}^{-i}\xi)\\
& = & \sum_{h=1}^{q-1}\sum_{i=1}^{\infty}\sum_{k\in\N_0}\hat{\psi^l}(\mathfrak{p}^{-j-1}(\mathfrak{p}\xi+u(k)))\overline{\hat{\psi}^h(\mathfrak{p}^{-i}(\mathfrak{p}\xi+u(k)))}\hat{\psi}^h(\mathfrak{p}^{-i}\mathfrak{p}\xi)\\
& = &  G^l_{j+1}(\mathfrak{p}\xi).
\end{eqnarray*}
This is equivalent to $G^l_j(\xi) = G^l_{j-1}(\mathfrak{p}^{-1}\xi)$. Iterating this equation, we obtain, $G^l_j(\xi) = G^l_1(\mathfrak{p}^{-j+1}\xi)$. We now calculate $G_1^l(\xi)$. We have
\begin{eqnarray*}
G_1^l(\xi)
& = & \sum_{k\in\N_0}\hat{\psi^l}(\mathfrak{p}^{-1}(\xi+u(k)))\sum_{h=1}^{q-1}\sum_{i=1}^{\infty}\overline{\hat{\psi}^h(\mathfrak{p}^{-i}(\xi+u(k)))}\hat{\psi}^h(\mathfrak{p}^{-i}\xi)\\
& = & \sum_{k\in \N_0}\hat{\psi^l}(\mathfrak{p}^{-1}\xi+u(qk))\sum_{h=1}^{q-1}\sum_{i=0}^{\infty}\overline{\hat{\psi}^h(\mathfrak{p}^{-i}(\mathfrak{p}^{-1}\xi+u(qk)))}\hat{\psi}^h(\mathfrak{p}^{-i}\mathfrak{p}^{-1}\xi) \\
& = & \sum_{k\in q\N_0}\hat{\psi^l}(\mathfrak{p}^{-1}\xi+u(k))\sum_{h=1}^{q-1}\sum_{i=0}^{\infty}\overline{\hat{\psi}^h(\mathfrak{p}^{-i}(\mathfrak{p}^{-1}\xi+u(k)))}\hat{\psi}^h(\mathfrak{p}^{-i}\mathfrak{p}^{-1}\xi) \\
\end{eqnarray*}
In the last equation, we can replace the sum over $k\in q\N_0$ by a sum over $k\in\N_0$ since the additional terms corresponding to $k\in\N_0\setminus q\N_0$ are all zero by (\ref{f5}). Hence,
\begin{eqnarray*}
 G^l_1(\xi)
& = & \sum_{h=1}^{q-1}\sum_{i=0}^{\infty}\hat{\psi}^h(\mathfrak{p}^{-i}\mathfrak{p}^{-1}\xi) \delta_{i,0}\delta_{l,h}\\
& = & \hat{\psi}^l(\mathfrak{p}^{-1}\xi).
\end{eqnarray*}
Thus $G^l_j(\xi) = \hat{\psi^l}(\mathfrak{p}^{-j}\xi)$ a.e. $\xi\in K$, since $\langle \omega_j^l(\xi), \omega_i^h(\xi)\rangle$ is integral-periodic, (\ref{omega}) follows. This completes the proof of Lemma.
\qed

We will also need the following lemma. We refer to~\cite{A}, \cite{CG} and~\cite{HW} for a proof of this result.
\begin{lemma}\label{lem:dim}
Let $\{v_j:j\geq 1\}$ be a family of vectors in a Hilbert space $\mathbb{H}$ such that
\begin{itemize}
\item [(i)] $\sum_{n=1}^{\infty}\|v_n\|^2 = C < \infty$, and
\item [(ii)] $v_n = \sum_{m=1}^{\infty}\langle v_n, v_m\rangle v_m$ for all  $n\geq 1$.
\end{itemize}
Let $\mathbb{F}=\overline{\rm span}\{v_j:j\geq 1\}$. Then $\dim\mathbb{F}=\sum_{j=1}^{\infty}\|v_j\|^2=C$.
\end{lemma}

The main result of this section is the following theorem.
\begin{theorem}
A wavelet $\Psi=\{\psi^1,\psi^2,\dots,\psi^{q-1}\}\subseteq L^2(K)$ is an MRA wavelet if only if $D_{\Psi}(\xi) = 1$ for almost every $\xi\in\ K$.
\end{theorem}

\proof
We have already observed that $D_{\Psi}(\xi) = 1$ for a.e. $\xi\in K$ when $\Psi$ is an MRA-wavelet. We now prove the converse. 

Assume that $D_{\Psi}(\xi) = 1$ for a.e. $\xi\in K$. Let $E$ be the subset of $\mathfrak D$ on which $D_{\Psi}(\xi)$ is finite and (\ref{omega}) is satisfied. Then $\omega_j^l$ are well-defined on $E$. For $\xi\in E$, we define the space
\[
\mathcal{F}(\xi)=\overline{\rm span}\{\omega^l_j(\xi): 1\leq l\leq q-1,j\geq 1\}.
\]
Then, by Lemma~\ref{lem:rel} and \ref{lem:dim}, we have
\begin{eqnarray}\label{f13}
\dim\mathcal{F}(\xi) = \sum_{l=1}^{q-1}\sum_{j=1}^{\infty}\|\omega_j^l(\xi)\|^2_{\ell^2}=D_{\Psi}(\xi)=1.
\end{eqnarray}
That is, for each $\xi\in E$, $\mathcal{F}(\xi)$ is generated by a single unit vector $U(\xi)$. We now choose a suitable vector. For $j\geq 1$, let us define
\[
X_j = \{\xi\in E: \omega_j^l(\xi)\neq 0 ~\mbox{for}~\mbox{some}~l,1\leq l\leq q-1,~\mbox{and}~\omega_m^l(\xi) = 0, \forall m<j ~\mbox{and}~ 1\leq l\leq q-1\},
\]
and 
\[
X_0 = \{\xi\in\mathfrak{D}:\omega_j^l(\xi)= 0 ~\mbox{for}~ ~\mbox{all}~ l, 1\leq l\leq q-1,~\mbox{and}~~\mbox{for}~ ~\mbox{all}~j\geq 1\}.
\]
Then $\{X_j:j=0, 1,2,\dots\}$ forms a partition of $E$. Note that $X_0=\{\xi\in\mathfrak{D}: D_\Psi(\xi)=0\}$. So for a.e. $\xi\in E\setminus X_0$, there exists $j\geq 1$ such that $\xi\in X_j$. Hence, there exists at least one $l$, $1\leq l\leq q-1$, such that $\omega^l_j(\xi)\neq 0$. Choose the smallest such $l$ and define
\[
U(\xi) = \frac{\omega_j^l(\xi)}{\|\omega_j^l(\xi)\|_{\ell^2}}.
\]
Thus, $U(\xi)$ is well defined and $\|U(\xi)\|_{\ell^2}=1$ for a.e. $\xi\in\mathfrak{D}$. We write $ U(\xi) = \{u_k(\xi):k\in\N_0\}$. Now, define $\hat{\varphi}(\xi) = u_k(\xi-u(k))$, where $k$ is the unique integer in $\N_0$ such that $\xi\in{\mathfrak D}+u(k)$. This defines $\hat\varphi$ on $K$. We first show that $\varphi\in L^2(K)$ and $\{\varphi(\cdot-u(k):k\in\N_0)\}$ is an orthonormal system in $L^2(K)$. We have
\begin{eqnarray*}
\|\hat{\varphi}\|^2_2
& = & \int_K |\hat{\varphi}(\xi)|^2d\xi \\
& = & \int_{\mathfrak{D}}\sum_{k\in\N_0}|\hat{\varphi}(\xi+u(k))|^2d\xi \\
& = & \sum_{k\in\N_0}\int_{\mathfrak{D}}|u_k(\xi)|^2 d\xi \\
& = & \int_{\mathfrak{D}}\|U(\xi)\|^2_{\ell^2}d\xi\\
& = & 1.
\end{eqnarray*}
Thus, $\varphi\in L^2(K)$. Also, 
\begin{equation}\label{f15}
\sum_{k\in\N_0} |\hat{\varphi}(\xi+u(k))|^2 = \sum_{k\in\N_0}|u_k(\xi)|^2 = \|U(\xi)\|^2_{\ell^2} = 1.
\end{equation}
This is equivalent to the fact that $\{\varphi(\cdot-u(k):k\in\N_0)\}$ is an orthonormal system. We now define $V_0^\sharp = \overline{\rm span}\{\varphi(\cdot-u(k)):k\in\N_0\}$. Let $W_j=\overline{\rm span}\{\psi^l_{j,k}:1\leq l\leq q-1, k\in\N_0\}$ and $V_0=\bigoplus\limits_{j<0}W_j$. If we can show that $V_0^\sharp=V_0$, then it will follow that $\{V_j: j\in\Z\}$ is the required MRA (see the discussion just after the definition of MRA). 

We first show that $V_0\subset V_0^\sharp$. It is sufficient to verify that $\psi^l_{j,k}\in V_0^\sharp$, $k\in\N_0$, $j<0$, $1\leq l\leq q-1$. For each $j\geq 1$, there exists a measurable function $\nu_j^l$ on $\mathfrak{D}$ such that $\omega_j^l(\xi) = \nu_j^l(\xi)U(\xi)$ for a.e. $\xi\in\mathfrak{D}$. That is,
\[
\hat{\psi}^l(\mathfrak{p}^{-j}(\xi+u(k))) = \nu_j^l(\xi)\hat{\varphi}(\xi+u(k))\quad\mbox{for a.e.}~\xi\in{\mathfrak D}, k\in\Z.
\]
Hence, by (\ref{f15}), for a.e. $\xi\in\mathfrak{D}$, we have
\begin{equation}\label{f16}
\sum_{k\in\N_0}|\hat{\psi}^l(\mathfrak{p}^{-j}(\xi+u(k)))|^2 = \sum_{k\in\N_0}|\nu_j^l(\xi)|^2|\hat{\varphi}(\xi+u(k))|^2 = |\nu_j^l(\xi)|^2.
\end{equation}
This shows that $\nu_j^l\in L^2(\mathfrak{D})$ so that we can write its Fourier series expansion. Thus, for $j\geq 1$, there exists $\{a_{j,k}^l:k\in\N_0\}\in\ell^2(\N_0)$ such that $\nu_j^l(\xi)=\sum_{k\in\N_0}a_{j,k}^l\overline{\chi_k(\xi)}$, with convergence in $L^2(\mathfrak{D})$. Extending $\nu_j^l$ integer periodically, we have
\begin{equation}\label{f17}
\hat{\psi}^l(\mathfrak{p}^{-j}\xi) = \nu_j^l(\xi)\hat{\varphi}(\xi) ~\mbox{for~a.e.}~ \xi\in K, j\geq 1.
\end{equation}
Taking inverse Fourier transform, we get
\[
\psi_{-j,0}^l(x) = q^{j/2}\sum_{k\in\N_0}a_{j,k}^l\varphi(x-u(k)), \quad j\geq 1.
\]
Hence, $\psi_{-j,0}^l\in V_0^\sharp$ for $j\geq 1$. Moreover, since $V_0^\sharp$ is invariant under translations by $u(k)$, $k\in\N_0$, we have $\psi^l_{j,k}\in V_0^\sharp$, $j<0$, $k\in\N_0$, $1\leq l\leq q-1$.

To show the reverse inclusion, it suffices to show that $V_0^\sharp\perp W_j$, for $j\geq 0$. For $j\geq 0, k\in\N_0, 1\leq l\leq q-1$, we have
\begin{eqnarray}\label{e.last}
\langle \varphi, \psi^l_{j,k}\rangle = \langle \hat\varphi, (\psi^l_{j,k})^\wedge\rangle
& = & \int_K\hat{\varphi}(\xi)q^{-j/2}\overline{\hat{\psi}^l(\mathfrak{p}^j\xi)}\chi_k(\mathfrak{p}^j\xi)d\xi \nonumber\\
& = & q^{j/2}\int_K\hat{\varphi}(\mathfrak{p}^{-j}\xi)\overline{\hat{\psi}^l(\xi)}\chi_k(\xi)d\xi \nonumber\\
& = & q^{j/2}\int_{\mathfrak{D}}\sum_{n\in\N_0}\hat{\varphi}(\mathfrak{p}^{-j}(\xi+u(n)))\overline{\hat{\psi}^l(\xi+u(n))}\chi_k(\xi)d\xi. 
\end{eqnarray}
Using equation (\ref{f16}), we get
\[
\sum_{l=1}^{q-1}\sum_{j=1}^{\infty}|\nu_j^l(\xi)|^2 = \sum_{l=1}^{q-1}\sum_{j=1}^{\infty}\sum_{k\in\N_0}|\hat{\psi}^l(\mathfrak{p}^{-j}(\xi+u(k)))|^2 = 1~\mbox{for~a.e.}~ \xi\in K.
\]
Hence, for such $\xi$ and for all $j\geq0$, there exists $j_0\geq 1$ such that $\nu^l_{j_0}(\mathfrak{p}^{-j}\xi)\neq 0$. Thus, (\ref{f17}) implies that $\hat{\psi}^l(\mathfrak{p}^{-j-j_0}\xi)=\nu^l_{j_0}(\mathfrak{p}^{-j}\xi)\hat{\varphi}(\mathfrak{p}^{-j}\xi)$. Therefore, for $k\in\N_0$, we get
\[
\hat{\psi}^l(\mathfrak{p}^{-j-j_0}(\xi+u(k)))
=\nu^l_{j_0}(\mathfrak{p}^{-j}(\xi+u(k)))\hat{\varphi}(\mathfrak{p}^{-j}(\xi+u(k))).
\]
Since $\mathfrak{p}^{-j}(\xi+u(k))=\mathfrak{p}^{-j}\xi+u(q^jk)$ and $\nu^l_{j_0}$ is integral-periodic, we have
\[
\hat{\varphi}(\mathfrak{p}^{-j}(\xi+u(k))) 
=\frac{1}{\nu^l_{j_0}(\mathfrak{p}^{-j}\xi)}\hat{\psi}^l(\mathfrak{p}^{-j-j_0}(\xi+u(k))).
\]
Therefore, using Lemma~\ref{lem:org}, for any $h$ with $1\leq h\leq q-1$, we have
\begin{eqnarray*}
\sum_{k\in\N_0}\hat{\varphi}(\mathfrak{p}^{-j}(\xi+u(k)))\overline{\hat{\psi}^h(\xi+u(k))}
& = & \frac{1}{\nu^l_{j_0}(\mathfrak{p}^{-j}\xi)}\sum_{k\in\N_0}\hat{\psi}^l(\mathfrak{p}^{-j-j_0}(\xi+u(k)))\overline{\hat{\psi}^h(\xi+u(k))}\\
& = & 0,
\end{eqnarray*}
since $j+j_0\geq 1$. Substituting this in~(\ref{e.last}), we get $\langle\varphi, \psi^l_{j,k}\rangle = 0$ for $j\geq 0, k\in\N_0, 1\leq l\leq q-1$. From this we conclude that $V_0^\sharp\subset V_0$. This completes the proof of Theorem.
\qed


\end{document}